\documentclass[preprint,review,number]{elsarticle}

\usepackage[utf8x]{inputenc}
\usepackage{cmap}
\usepackage[T1]{fontenc}
\usepackage[british]{babel}
\usepackage[english=british]{csquotes}
\usepackage{hyphenat}

\usepackage{units}
\usepackage{amsmath} 
\usepackage{amssymb}
\usepackage{amsthm}
\usepackage{stmaryrd}
\usepackage{mathtools}

\usepackage{booktabs}
\usepackage{caption}

\usepackage{graphicx}
\graphicspath{{./pictures/},{./}}
\DeclareGraphicsExtensions{.pdf}

\newcommand{\ddx}[2]{\frac{{\mathrm d}{#1}}{{\mathrm d}{#2}}}
\newcommand{\DDx}[2]{\frac{{\mathrm D}{#1}}{{\mathrm D}{#2}}}
\newcommand{\pdx}[2]{\frac{\partial {#1}}{\partial {#2}}}
\newcommand{\idx}[4]{\int_{#1}^{#2} {#3}\,{\mathrm d}{#4}}

\newcommand{\vv}[1]{\mathbf{#1}}        
\newcommand{\gv}[1]{\boldsymbol{#1}}    
\newcommand{\sm}[1]{\mathsf{#1}}        
\newcommand{\sv}[1]{\mathsf{#1}}        

\newcommand{\bigO}[1]{\mathcal{O}(#1)}
\newcommand{\qnorm}[1]{\left\Vert{#1}\right\Vert_{W^{-k,q}}}
\newcommand{\pnorm}[1]{\left\Vert{#1}\right\Vert_{W^{k,p}}}
\newcommand{\cnorm}[1]{\left\Vert{#1}\right\Vert_{C^{n+2}}}
\newcommand{\lnorm}[1]{\left\Vert{#1}\right\Vert_{l^{1}}}
\newcommand{\pspace}{W^{k,p}}

\newcommand{\pairing}[2]{\left\langle{#1},{#2}\right\rangle}
\newcommand{\Cemb}{C_{\mathrm{emb}}}
\newcommand{\neighi}{\mathcal{N}_i}

\newtheorem{lem}{Lemma}
\newtheorem{thm}{Theorem}
\newproof{pf}{Proof}

\begin{document}
%
\title{A Splitting-free Vorticity Redistribution Method}
\author{M.~Kirchhart\fnref{cor1}}
\ead{kirchhart@keio.jp}
\fntext[cor1]{Corresponding author}

\author{S.~Obi}
\ead{obsn@mech.keio.ac.jp}

\address{Department of Mechanical Engineering, Keio University,
3-14-1 Hiyoshi, Kōhoku-ku, Yokohama-shi 233-8522, Japan}

\begin{abstract}
We present a splitting-free variant of the vorticity redistribution method.
Spatial consistency and stability when combined with a time-stepping scheme are
proven. We propose a new strategy preventing excessive growth in the number of
particles while retaining the order of consistency. The novel concept of small
neighbourhoods significantly reduces the method's computational cost. In
numerical experiments the method showed second order convergence, one order
higher than predicted by the analysis. Compared to the fast multipole code
used in the velocity computation, the method is about three times faster.
\end{abstract}
\maketitle
%
%
%
\section{Introduction}
In vortex methods in two dimensions and in the absence of boundaries, one wants
to evolve a scalar vorticity field $\omega$ in form of a \emph{particle cloud}:
\begin{equation}\label{eqn:vfield}
\omega(t,\vv{x}) = \sum_{i=1}^{N}\Gamma_i(t)\delta(\vv{x}-\vv{x}_i(t)),
\end{equation}
over time $t$ according to the vorticity transport equation:
\begin{equation}\label{eqn:vte}
\DDx{\omega}{t}\equiv\pdx{\omega}{t} + (\vv{u}\cdot\nabla)\omega = \nu\Delta\omega.
\end{equation}
In there, $\Gamma_i\in\mathbb{R}$ denotes the circulation that particle $i$
carries, $\vv{x}_i\in\mathbb{R}^2$ stands for that particle's position, 
$\vv{u}:\mathbb{R}^2\to\mathbb{R}^2$ is the velocity field induced by $\omega$
according to the Biot--Savart law, $\nu\geq 0$ refers to the fluid's kinematic
viscosity, and $\delta$ is the Dirac delta distribution.

The beauty of vortex methods lies in their handling of the inviscid case
($\nu~=~0$): evolving $\Gamma_i$ and $\vv{x}_i$ according to the following set
of ordinary differential equations (ODEs):
\begin{equation}\label{eqn:eulerode}
\begin{split}
\ddx{\vv{x}_i}{t} &= \vv{u}(t,\vv{x}_i(t)), \\
\ddx{\Gamma_i}{t} &= 0,
\end{split}
\end{equation} 
i.\,e., by convecting the particles according to the local velocity and leaving
their strengths unchanged, the resulting vorticity field fulfils the vorticity
transport equation~\eqref{eqn:vte} exactly. Especially, due to the absence of a
fixed computational mesh and the natural treatment of convection, inviscid vortex
methods are \emph{free of numerical dissipation} and conserve circulation, linear
and angular momentum, as well as energy exactly~\cite{cottet2000}. 

Many different approaches on how to handle the viscous case have been suggested
in the literature, the book by Cottet and Koumoutsakos~\cite{cottet2000} gives
an overview and references to some of the most commonly used approaches. Almost
all of them belong to the class of \emph{viscous-splitting} algorithms: first,
particles are convected under the absence of viscosity. Afterwards vorticity
is diffused according to the heat-equation, i.\,e., in absence of convection.
One of the earliest such approaches is the so-called \enquote{random-walk}
method, in which viscosity is simulated by an additional Brownian motion of
the particles. This method, however, converges only very slowly.

The resurrected core-spreading technique~\cite{rossi1996} relies on a different
representation of the vorticity field~\eqref{eqn:vfield}. The Biot--Savart kernel
is singular at the origin, causing very large velocity values when particles
approach each other. For this reason the Biot--Savart law is usually regularised
in practice. This is commonly done by replacing the Dirac delta distribution
$\delta$ with a smooth approximation $\zeta_\varepsilon$, a so-called blob-%
function with blob- or core-width $\varepsilon$. As the name suggests, the core-%
spreading method works by enlarging individual particle's core widths~%
$\varepsilon$. This enlargement causes the solution to get increasingly blurred
over time, unless some kind of remeshing is employed.

The method of particle strength exchange~(PSE), on the other hand, modifies the
particle strengths by approximating the Laplacian by an integral. This integral
is then approximated by numerical quadrature, using the particle positions as
quadrature nodes. Frequent remeshing is required, unless the newer mesh-free
variant DC-PSE~\cite{schrader2010} is employed.

The vorticity redistribution method (VRM) by Shankar and van~Dommelen~%
\cite{shankar1996} can be interpreted as a computed finite-difference
stencil which solves the heat-equation for a given time-step~$\Delta t$. The
fact that these stencils are computed on-the-fly makes the method completely
mesh-free. 

While most of these methods achieve high-order spatial accuracy, the viscous\hyp%
splitting inevitably limits their accuracy in time to first order, unless more
sophisticated splits are used~\cite{beale1981}. Note that this result holds
regardless of the time-stepping scheme used, underlining that splitting the
equation is unnatural: diffusion and convection \emph{do happen simultaneously}
and thus should not be treated one after another.

Our contribution in this article is a new method of treating the diffusive term
in a manner similar to the vorticity redistribution method (VRM). Instead of
computing a stencil that approximates a solution to the heat equation, we
directly approximate the Laplacian. This allows us to avoid the viscous-%
splitting and to treat both diffusion and convection simultaneously. The spatial
consistency of our method is proven. As we do not need to integrate the Laplacian
over time, our proof does not require Fourier analysis, like the original VRM.
We then consider the case of pure diffusion in combination with the forward
Euler method and derive sharp a-priori and a-posteriori bounds on the
step-width. This analysis in the absence of convection is justified, as the
convective part of the equations is known to be stable independent of the
step-width~\cite{anderson1985}. The resulting a-priori bound is---apart from a
constant---identical to the classical stability condition for the five-point
central-difference stencil, underlining the interpretation of our method as a
computed finite-difference method.

Finally, we show that our method conserves circulation, linear, and angular
momentum. In the original VRM-paper it was suggested to ignore particles in the
diffusive process if their circulation was below a certain threshold. Choosing
a low threshold does yield accurate discretisations, however, the choice of its
value seemed rather arbitrary. We propose a new strategy preventing excessive
growth in the number of particles while maintaining the order of consistency.
Based on results by Seibold~\cite{seibold2006,seibold2008}, we further introduce
the new concept of \emph{small neighbourhoods} which significantly reduces the
computational cost of the method. The resulting scheme keeps all of the benefits of
the original VRM while not relying on viscous splitting or arbitrary thresholds.
We conclude with numerical examples illustrating efficiency and convergence of
the method in the purely diffusive, as well as in the convective case.

\section{Description of the Method}
Our aim is to approximate the Laplacian of $\omega$ by the following formula:
\begin{equation}
\Delta_h\omega := \sum_{i=1}^{N}\sum_{j=1}^{N}f_{ij}\Gamma_i\delta(\vv{x}-\vv{x}_j),
\end{equation}
where $f_{ij}\Gamma_i$ refers to the rate at which circulation is diffused from
particle $i$ to particle $j$. The values of $f_{ij}$ need to be chosen such that
certain conditions are fulfilled in order for this approximation to be accurate.

In order to specify these conditions, we assume that the particles are quasi-%
uniformly distributed, with $h$ corresponding to the average inter-particle
spacing. We then define the neighbourhood $\neighi$ of particle $i$ as follows:
\begin{equation}\label{eqn:neighbourhood}
\neighi := \lbrace j\in{1,\ldots,N}:
rh\leq |\vv{x}_i-\vv{x}_j|\leq Rh\rbrace\cup\lbrace i\rbrace.
\end{equation}
where $R>r>0$ are fixed, user-defined parameters. The original VRM formulation
does not include the lower bound $r$. Due to their movement, particles might get
closer to one another than this. In section~\ref{sec:outlook} we will give some
remarks on this problem. In our analysis we show that both bounds are required
to control the error: the upper bound limits the cut-off error of the expansions
used, while the lower bound is needed for stability. For $j\in\neighi$ the
values $f_{ij}$ are chosen such that certain moment conditions are fulfilled.
For $j\not\in\neighi$ we define $f_{ij} = 0$.

As will be seen later on, depending on the particle cloud's geometry, these
moment conditions do not always have a solution. Note, however, that we can
always add new particles of zero strength to the field: introducing a new
particle with $\Gamma_i = 0$ leaves the vorticity field~\eqref{eqn:vfield}
unchanged. In this point the VRM sligthly differs from classical finite-difference
methods: the vorticity field~\eqref{eqn:vfield} is not a list of pointwise
function values; rather it can be seen as a quadrature rule for integrating
functions against an underlying, smooth vorticity field. Inserting an empty
particle corresponds to adding a quadrature node with weight zero.

For such empty particles one obviously always has $\Gamma_i = 0\Longrightarrow
f_{ij}\Gamma_i = 0$, i.\,e., the value of $f_{ij}$ is arbitrary and can safely
be defined as zero, too. We make use of this fact by inserting new particles to
fill holes in the cloud and to expand it at its outer rim. This way we can
ensure that for all circulation-carrying particles sufficiently many neighbours
do exist. Circulation will then be diffused to the neighbouring particles and
thereby be spread out in space, which also is in accordance with the physical
intuition of diffusive processes.

At the core of our method lies the computation of the values $f_{ij}$ for every
$i$ and $j\in\neighi$. In order to ensure accuracy, the error is developed as a
Taylor expansion. We require that at least all error terms of constant, linear,
and quadratic order vanish. For second order accuracy one may also choose to
require cubic terms to vanish. A detailed derivation of the resulting equations
is given in section~\ref{sec:consistency}.

As will be seen later on, non-negativity of stencils is a sufficient criterion
for stable time-discretisations. In addition to that, such stencils posses many
more desirable properties, as described by Seibold~\cite{seibold2006,seibold2008,seibold2010}.
A stencil is called non-negative if it fulfils $f_{ij}\geq 0$ for all $j\neq i$.
Unfortunately, as also will be shown in the analysis section~\ref{subsec:limitations},
non-negative stencils cannot fulfil the moment equations of fourth order. Unless
one gives up on non-negativity and the resulting stability guarantee, the method's
accuracy is hence limited to second order.

Like the two-dimensional Taylor expansion, the moment conditions are most easily
expressed using multi-index notation. Defining the vector $\vv{r}_{ij}$:
\begin{equation}
\vv{r}_{ij} = \vv{x}_j - \vv{x}_i, 
\end{equation}
and denoting its Cartesian components by $\vv{r}_{ij}^x$ and $\vv{r}_{ij}^y$,
respectively, for $\bigO{h^n}$ accuracy, with $n=1$ or $n=2$, we pose the
following conditions:
\begin{equation}
\label{eqn:cond_two}
\sum_{j=1}^{N}f_{ij}\vv{r}_{ij}^x\vv{r}_{ij}^x = 2,\quad
\sum_{j=1}^{N}f_{ij}\vv{r}_{ij}^y\vv{r}_{ij}^y = 2,\quad 
\sum_{j=1}^{N}f_{ij}\vv{r}_{ij}^x\vv{r}_{ij}^y = 0,
\end{equation}
and for all other error terms with multi-index $\alpha$:
\begin{equation}
\label{eqn:cond_alpha}
\sum_{j=1}^{N}f_{ij}\vv{r}_{ij}^\alpha = 0, \qquad
0\leq |\alpha|\leq n + 1,\ |\alpha|\neq 2.
\end{equation}
Because we have $\vv{r}_{ii}\equiv\boldsymbol{0}$, only the equation for
$\alpha=(0,0)$ depends on $f_{ii}$, yielding:
\begin{equation}\label{eqn:fii}
f_{ii} = -\sum_{j\neq i}f_{ij}.
\end{equation}
For $n=1,2$ we consequently have to solve a system consisting of five or nine
moment conditions, respectively. For every particle $i$, this linear system can
be rewritten in matrix-vector notation:
\begin{equation}\label{eqn:momcondmv}
\sm{V}_i\sv{f}_i = \sv{b},\ \sv{f}_i \geq 0.
\end{equation}
Here, $\sv{f}_i$ is the vector of coefficients $f_{ij}$, $i\neq j$, $\sv{b}$ is
the vector that contains only zero entries except for the two~\enquote{2}-entries
at $\alpha=(2,0)$ and $\alpha=(0,2)$, and $\sm{V}_i$ is the Vandermonde-matrix,
with rows for each multi-index $1\leq|\alpha|\leq n + 1$ and columns $j$ for
each particle $j\in\neighi\setminus\lbrace i\rbrace$:
\begin{equation}
V_{\alpha,j} = \vv{r}_{ij}^\alpha.
\end{equation}
In order to obtain scaling independent of $h$, for a numerical implementation
it is beneficial and straightforward to rewrite these conditions for the
normalised vectors $\vv{r}_{ij}/h$. In section~\ref{sec:implementation} we
describe how to solve these equations and how to ensure that non-negative
stencils exist.

\section{Analysis}\label{sec:analysis}
\subsection{Preliminaries}
Let $n\in\lbrace1,2\rbrace$ be the desired order of accuracy, let $p\in[1,\infty)$
be arbitrary but fixed, let $q$ be its conjugate exponent such that $1=1/p + 1/q$,
and let $k$ be an integer such that $k>2/p + n + 2$. We denote the Sobolev
space of $k$ times weakly differentiable $L^p(\mathbb{R}^2)$-functions by
$W^{k,p}$, and let $W^{-k,q}$ refer to its dual space. Note that, by the Sobolev
embedding theorem, we have $W^{k,p}\hookrightarrow C^{n+2}$, where $C^{n+2}$
refers to the space of $n+2$ times continuously differentiable functions equipped
with the maximum norm over all derivatives. Further note that we have
$\Vert(\Gamma_i)\Vert_{l^1}<\infty$, and therefore $\omega(t,\cdot)\in W^{-k,q}$:
\begin{equation}\label{eqn:normbound}
\begin{split}
\qnorm{\omega} &= \sup_{\varphi\in\pspace}
\frac{\pairing{\omega}{\varphi}}{\pnorm{\varphi}}
\leq \sup_{\varphi\in\pspace}\Cemb
     \frac{\sum_{i=1}^{N}|\Gamma_i|\,|\varphi(\vv{x}_i)|}{\cnorm{\varphi}} \\
&\leq \sup_{\varphi\in\pspace}\Cemb \frac{\cnorm{\varphi}\sum_{i=1}^{N}|\Gamma_i|}{\cnorm{\varphi}}
 = \Cemb\Vert(\Gamma_i)\Vert_{l^1},
\end{split}
\end{equation}
where $\pairing{\cdot}{\cdot}$ refers to the dual pairing and $\Cemb$ denotes the
Sobolev embedding constant. This inequality also allows us to infer stability in
the $\qnorm{\cdot}$-norm by bounding the $l^1$-norm of the circulations later on.

\subsection{Consistency}\label{sec:consistency}
We will need the following lemma.
\begin{lem}\label{lem:fnorm}
For a stencil that satisfies~\eqref{eqn:momcondmv} one has:
\begin{equation}
f_{ii} = -\sum_{j\neq i} f_{ij} \leq 0,\quad 
4(Rh)^{-2} \leq \sum_{j\neq i}f_{ij} \leq 4(rh)^{-2}.
\end{equation}
\end{lem}
\begin{pf}
The first part directly follows from $f_{ij}\geq 0$ for $i\neq j$ and
equation~\eqref{eqn:fii}. The second relation follows by the sum of the moment
equations for $\alpha=(2,0)$ and $\alpha=(0,2)$ and $rh \leq |\vv{r}_{ij}|\leq Rh$.\qed
\end{pf}
We now are ready to prove the following consistency result.
\begin{thm}[Consistency]\label{thm:consistency}
One has:
\begin{equation*}
\qnorm{\Delta\omega - \Delta_h\omega} \leq
C\biggl(\frac{R}{r}\biggr)^2 (Rh)^{n}\lnorm{(\Gamma_i)},
\end{equation*}
where $C$ is a constant that only depends on $n$.
\end{thm}
\begin{pf}
For arbitrary $\varphi\in\pspace$ one has:
\begin{equation}
\label{eqn:error_start}
\pairing{\Delta\omega - \Delta_h\omega}{\varphi} = 
\sum_{i=1}^{N}\Gamma_i\biggl(\Delta\varphi(\vv{x}_i) - 
\sum_{j=1}^{N}f_{ij}\varphi(\vv{x}_j)\biggl).
\end{equation}
We develop $\varphi(\vv{x}_j)$ as a Taylor series around $\varphi(\vv{x}_i)$ and
obtain:
\begin{equation}
\varphi(\vv{x}_j) = \sum_{|\alpha|\leq n+1}\frac{\vv{r}_{ij}^\alpha}{\alpha!}D^{\alpha}\varphi(\vv{x}_i)
                  + \underbrace{
                    \sum_{|\alpha| =   n+2}\frac{\vv{r}_{ij}^\alpha}{\alpha!}D^{\alpha}\varphi(\gv{\xi}_{ij}),
                     }_{=:R_{ij}^{n+2}}               
\end{equation}
where $\gv{\xi}_{ij}$ is a point on the line connecting $\vv{x}_i$ and
$\vv{x}_j$. The moment conditions were chosen such that the first sum vanishes
when this relation is inserted into equation~\eqref{eqn:error_start}. Note
that we have with the help of the Sobolev embedding:
\begin{equation}
|R_{ij}^{n+2}|\leq \Cemb C_\alpha(Rh)^{n+2}\pnorm{\varphi},\ 
C_\alpha = \sum_{|\alpha|=n+2}\frac{1}{\alpha!},
\end{equation}
such that we get with the help of the triangle inequality, H\"older's inequality,
and Lemma~\ref{lem:fnorm}:
\begin{equation}
\begin{split}
|\pairing{\Delta\omega - \Delta_h\omega}{\varphi}| &=
\biggl\vert\sum_{i=1}^{N}\Gamma_i\sum_{j=1}^{N}f_{ij}R_{ij}^{n+2}\biggr\vert
\leq  \lnorm{(\Gamma_i)}\bigl\Vert(\sum_{j=1}^{N}f_{ij}R_{ij}^{n+2})_i\bigr\Vert_{l^\infty} \\
&\leq  4\Cemb C_\alpha \frac{(Rh)^{n+2}}{(rh)^2}\pnorm{\varphi}\lnorm{(\Gamma_i)}.
\end{split}
\end{equation}\qed
\end{pf}

\subsection{Stability for the Heat Equation}
In our next step we investigate the stability of Euler’s method in combination
with our spatial discretisation. As we introduced a new discretisation of the
Laplace operator, it is natural to omit convection and to investigate the heat
equation:
\begin{equation}\label{eqn:heat_eq}
\pdx{\omega}{t} = \nu\Delta\omega.
\end{equation}
We are going to apply the method of lines: in our case $\omega$ is a particle
cloud, the Laplacian operator is replaced with its discretisation $\Delta_h$,
and the time derivative is discretised using Euler's method. While this method
is only first-order accurate, this all that is needed to construct higher
order schemes: so-called non-linear SSP-stable methods of higher order exist,
which can be written as a convex combination of several Euler steps~\cite{gottlieb2001}.
While the classical Runge--Kutta method (RK4) is not such a scheme, our numerical
experiments exhibited no instabilities.

To ease notation, we introduce the vector $\sv{\Gamma}\in\mathbb{R}^N$,
consisting of the components $\Gamma_i$, and the matrix $\sm{F}\in\mathbb{R}^{%
N\times N}$, consisting of components $f_{ij}$, respectively. Denoting the
current and next time-steps with $n$ and $n+1$, respectively, our scheme then
reads:
\begin{equation}
\sv{\Gamma}^{n+1} = \underbrace{(\sm{I}+\nu\Delta t\sm{F}^\top)}_{=:\sm{C}}\sv{\Gamma}^n,
\end{equation}
where $\sm{I}\in\mathbb{R}^{N\times N}$ is the identity matrix and $\Delta t>0$
denotes the step-width. As shown in theorem~\ref{thm:consistency}, the consistency
error can be bounded by $\Vert\sv{\Gamma}\Vert_1$. It is therefore sufficient to
require $\Vert\sm{C}\Vert_1\leq 1$. Note that due to equation~\eqref{eqn:normbound},
this implies that $\qnorm{\omega}$ remains bounded as well. The following theorem
will show that positive stencils are not only sufficient but also necessary to
obtain a scheme that fulfils $\Vert\sm{C}\Vert_1\leq 1$. 

\begin{thm}[Stability]\label{thm:stability}
One has:
\begin{equation*}
\Vert\sm{C}\Vert_1 = 1,
\end{equation*}
if and only if we have a positive stencil:
\begin{equation*}
f_{ii}\leq 0,\quad f_{ij} \geq 0\ (i\neq j)
\end{equation*}
and for all $i=1,\ldots, N$:
\begin{equation*}
\nu\Delta t \leq -f_{ii}^{-1}.
\end{equation*}
For larger $\Delta t$ or non-positive stencils one always has
$\Vert\sm{C}\Vert_1>1$.
\end{thm}
\begin{pf}
One has:
\begin{equation}
\Vert\sm{C}\Vert_1 = \max_{j} \sum_{i=1}^{N}|C_{ij}| =
\max_{j} |1+\nu\Delta tf_{jj}| + \nu\Delta t\sum_{i\neq j}|f_{ji}|.
\end{equation}
Thus $\Vert\sm{C}\Vert_1\leq 1 \Longrightarrow f_{jj}\leq 0$. Now assume
$(1+\nu\Delta tf_{jj}) \geq 0$, i.\,e., $\nu\Delta t \leq -f_{jj}^{-1}$. We then
have for each $j$, due to equation~\eqref{eqn:fii}:
\begin{equation}
\sum_{i=1}^{N}|C_{ij}| =  1 - \nu\Delta t\sum_{i\neq j}f_{ji} +
                              \nu\Delta t\sum_{i\neq j}|f_{ji}|.
\end{equation}
Thus, we have $\Vert\sm{C}\Vert_1\leq 1$ if and only if for all $j$:
\begin{equation}
\sum_{i\neq j}|f_{ji}|\leq \sum_{i\neq j}f_{ji} \iff f_{ji}\geq 0.
\end{equation}
For positive stencils both sides are equal, and thus $\Vert\sm{C}\Vert_1 = 1$.

Conversely assume $(1+\nu\Delta tf_{jj}) < 0$, i.\,e., 
$\nu\Delta t > -f_{jj}^{-1}$. We then have again due to equation~\eqref{eqn:fii}:
\begin{equation}
\sum_{i=1}^{N}|C_{ij}| = -1 + \nu\Delta t\sum_{i\neq j}(f_{ij}+|f_{ij}|).
\end{equation}
Assume we would have $\Vert\sm{C}\Vert_1\leq 1$. We then would have for all $j$:
\begin{equation}
\sum_{i\neq j}(f_{ij}+|f_{ij}|)\leq\frac{2}{\nu\Delta t}.
\end{equation}
But note that we have:
\begin{equation}
\sum_{i\neq j}(f_{ij}+|f_{ij}|) \geq 
2\sum_{i\neq j}f_{ij} = 2f_{jj},
\end{equation}
and thus:
\begin{equation}
2f_{jj}\leq \frac{2}{\nu\Delta t} \iff \nu\Delta t\leq -f_{jj}^{-1},
\end{equation}
which is a direct contradiction to our assumption on the time-step.\qed
\end{pf}
Theorem~\ref{thm:stability} gives us an easy a-posteriori bound which can readily
be implemented. This allows us to optimally choose the step-width in a computer
program. In higher-order Runge--Kutta schemes it is hard to predict the values
$f_{ii}$ for intermediate stages. Thus, again employing Lemma~\ref{lem:fnorm},
the following a-priori bound is useful:
\begin{equation}\label{eqn:apriobound}
\Delta t\leq \frac{(rh)^{2}}{4\nu}.
\end{equation}
Note that this closely resembles the classical stability condition for the
five-point finite-difference stencil, highlighting the similarity between the
two methods. The fact that we can only achieve $\Vert\sm{C}\Vert_1 = 1$, as
opposed to $\Vert\sm{C}\Vert_1 < 1$, can be seen as a consequence of the fact
that our method conserves circulation, as will be shown in the next section.
 
\subsection{Conservation Properties for the Navier--Stokes Equations}
We now discuss the conservation properties of our method when used in combination
with convection, i.\,e., for the Navier--Stokes equations. We thus consider the
following semi-discrete system of coupled ordinary differential equations:
\begin{equation}\label{eqn:nsode}
\begin{split}
\ddx{\vv{x}_i}{t} &= \sum_{j=1}^{N}\vv{K}(\vv{r}_{ij})\Gamma_j, \\
\ddx{\Gamma_i}{t} &= \nu\sum_{j=1}^{N}f_{ji}\Gamma_j,
\end{split}
\end{equation}
where $\vv{K}$ denotes the Biot--Savart kernel defined for $\vv{x} = (x,y)\in
\mathbb{R}^2$ as:
\begin{equation}
\vv{K}(\vv{x}) :=
\begin{cases}
\gv{0} & \text{if $\vv{x} = \gv{0}$,} \\
\frac{(y,-x)^\top}{2\pi|\vv{x}|^2} & \text{else.}
\end{cases}
\end{equation}
Note that we differ from the usual sign convention, such that we do not need to
negate the vector~$\vv{r}_{ij}$ in equation~\eqref{eqn:nsode}. Also note that the
fractions $f_{ij}$ depend on the particle positions, which we, for reasons of
brevity, did not introduce into the notation.
We are going to investigate the following quantities:
\begin{itemize}
\item \leavevmode\rlap{Circulation:}\phantom{Angular Momentum:}
$I_0 := \idx{\mathbb{R}^2}{}{\omega\ \ \;\,}{\vv{x}}   = \sum_{i=0}^{N}\Gamma_i$,
\item \leavevmode\rlap{Linear Momentum:}\phantom{Angular Momentum:}
$\vv{I}_1 := \idx{\mathbb{R}^2}{}{\omega\vv{x}\ \,}{\vv{x}} = \sum_{i=0}^{N}\Gamma_i\vv{x}_i$,
\item Angular Momentum:
$I_2 := \idx{\mathbb{R}^2}{}{\omega\vv{x}^2}{\vv{x}}   = \sum_{i=0}^{N}\Gamma_i\vv{x}_i^2$.
\end{itemize}
The conservation laws for these quantities read~\cite{majda2001,tao2014}:
\begin{equation}\label{eqn:conservationlaws}
\ddx{I_0}{t} = 0,\ \ddx{\vv{I}_1}{t} = \gv{0},\ \ddx{I_2}{t} = 4\nu I_0.
\end{equation}
Note that these quantities are moments of vorticity and thus are closely linked
to the moment conditions~\eqref{eqn:cond_alpha} and~\eqref{eqn:cond_two}. This
close link will allow us to show that the semi-discrete equations~%
\eqref{eqn:nsode} fulfil the conservation laws~\eqref{eqn:conservationlaws}
exactly.
\begin{thm}[Conservation of Circulation and Momentum]
The vorticity field described by the system of ODEs~\eqref{eqn:nsode} conserves
circulation as well as linear and angular momentum.
\end{thm}
\begin{pf}
Our proof utilises the moment conditions as well as ideas from Cottet and
Koumoutsakos~\cite{cottet2000}. For circulation we immediately obtain:
\begin{equation}
\ddx{I_0}{t} = \sum_{i=1}^{N}\ddx{\Gamma_i}{t} 
= \nu\sum_{i=1}^{N}\sum_{j=1}^{N}f_{ji}\Gamma_j 
= \nu\sum_{j=1}^{N}\Gamma_j \underbrace{\sum_{i=1}^{N}f_{ji}}_{=0} = 0.
\end{equation}
For linear momentum we have:
\begin{equation}
\ddx{\vv{I}_1}{t} =
\sum_{i=1}^{N}\Gamma_i\ddx{\vv{x}_i}{t} +
\sum_{i=1}^{N}\vv{x}_i\ddx{\Gamma_i}{t}.
\end{equation}
For the first part of the sum note that the Biot--Savart kernel is odd, i.\,e.,
we have $\vv{K}(\vv{r}_{ij}) = -\vv{K}(\vv{r}_{ji})$. Using this relation and
exchanging the indices we obtain:
\begin{equation}
\sum_{i=1}^{N}\Gamma_i\ddx{\vv{x}_i}{t} =
 \sum_{i=1}^{N}\sum_{j=1}^{N}\vv{K}(\vv{r}_{ij})\Gamma_i\Gamma_j =
-\sum_{i=1}^{N}\sum_{j=1}^{N}\vv{K}(\vv{r}_{ij})\Gamma_i\Gamma_j.
\end{equation}
Thus, this part of the sum equals its negative and therefore is zero. For the
second part we have using the moment conditions:
\begin{equation}
\nu\sum_{i=1}^{N}\sum_{j=1}^{N}f_{ji}\Gamma_j\vv{x}_i =
\nu\sum_{j=1}^{N}\Gamma_j
\biggl(%
\underbrace{\sum_{i=1}^{N}f_{ji}\vv{r}_{ji}}_{=\gv{0}} +
\vv{x}_j\underbrace{\sum_{i=1}^{N}f_{ji}}_{=0}%
\biggr)
= \gv{0}.
\end{equation}
Lastly, for the angular momentum we obtain:
\begin{equation}
\ddx{I_2}{t} = 
\sum_{i=1}^{N}2\Gamma_i\vv{x}_i\cdot\ddx{\vv{x}_i}{t} +
\sum_{i=1}^{N}\vv{x}_i^2\ddx{\Gamma_i}{t}.
\end{equation}
For the first sum we have:
\begin{equation}
2\sum_{i=1}^{N}\sum_{j=1}^{N}\Gamma_i\Gamma_j\vv{K}(\vv{r}_{ij})\cdot\vv{x}_i.
\end{equation}
By writing $\vv{x}_i = \nicefrac{1}{2}(\vv{x}_i+\vv{x}_j) + 
\nicefrac{1}{2}(\vv{x}_i-\vv{x}_j)$ this sum again splits up into two parts.
Using the oddness property of $\vv{K}$ and exchanging the indices as above, the
first part is zero. For the second part note that by the definition of $\vv{K}$
we have $\vv{K}(\vv{r}_{ij})\cdot\vv{r}_{ij}\equiv 0.$ Finally, we have
$\vv{x}_i^2 = \vv{r}_{ji}^2 + 2\vv{x}_i^x\vv{x}_j^x + 2\vv{x}_i^y\vv{x}_j^y -
\vv{x}_j^2$ and thus:
\begin{multline}
\sum_{i=1}^{N}\vv{x}_i^2\ddx{\Gamma_i}{t} =
\nu\sum_{i=1}^{N}\sum_{j=1}^{N}f_{ji}\Gamma_j\vv{x}_i^2 = \\
\nu\sum_{j=1}^{N}\Gamma_j
\biggl(%
\underbrace{\sum_{i=1}^{N}f_{ji}\vv{r}_{ji}^2}_{=4} +
2\vv{x}_j^x\underbrace{\sum_{i=1}^{N}f_{ji}\vv{x}_{i}^x}_{=0} +
2\vv{x}_j^y\underbrace{\sum_{i=1}^{N}f_{ji}\vv{x}_{i}^y}_{=0} -
\vv{x}_j^2\underbrace{\sum_{i=1}^{N}f_{ji}}_{=0}%
\biggr) = \\
4\nu\sum_{i=1}^{N}\Gamma_i = 4\nu I_0.
\end{multline}\qed
\end{pf}
Due to the non-linear coupling of $\Gamma_{i}$ and $\vv{x}_i$ in $\vv{I}_1$ and
$I_2$, these quantities are generally not exactly conserved when the system of
ODEs~\eqref{eqn:nsode} is discretised using Euler's method. Here, one can only
verify $I_0$ to be conserved exactly. The numerical experiments at the end of this
article have shown, however, that the other two quantities are conserved very
well in practice. 
 
At the end we want to remark that kinetic energy $E$:
\begin{equation}
E = \idx{\mathbb{R}^2}{}{\omega(y\vv{u}^{x}-x\vv{u}^{y})}{\vv{x}}
\end{equation}
contains a product with the velocity $\vv{u}$, which is non-linearly linked to
$\omega$  and the particle positions through the Biot--Savart law. As the moment
conditions do not reflect this, we cannot expect the scheme to conserve energy
exactly.

\subsection{Reducing the Number of Diffused Particles}\label{subsec:reduced}
When used as introduced above, the method may give rise to large numbers of
particles carrying negligible amount of circulation, thus unnecessarily
increasing the numerical cost. For this reason, Shankar and van Dommelen~%
\cite{shankar1996} suggest to only diffuse particles carrying more circulation
than a prescribed threshold. In their work, they set this threshold to the
machine epsilon for single-precision floating-point arithmetic, i.\,e., round
about $5.96\cdot10^{-8}$.

While choosing a threshold near machine accuracy does produce accurate results,
this choice remains rather arbitrary. It is also not clear how big the
introduced error is. Luckily our analysis of the error may be extended to allow
for the exclusion of particles from diffusion. Let $\mathcal{I}$ be the set of
particles that are not diffused and let $\tilde{\Delta}_h$ refer to the
corresponding \enquote{reduced} approximation of the Laplacian:
\begin{equation}\label{eqn:reducedop}
\tilde{\Delta}_h\omega(\vv{x}) :=
\sum_{i\not\in\mathcal{I}}\sum_{j=1}^{N}\Gamma_i f_{ij}\delta(\vv{x}-\vv{x}_j).
\end{equation}
Using the same techniques as above, we can then investigate the additional error
introduced:
\begin{thm}
One has:
\begin{equation*}
\qnorm{(\Delta_h-\tilde{\Delta}_h)\omega}\leq
4(rh)^{-2}\Cemb\lnorm{(\Gamma_i)_{i\in\mathcal{I}}}.
\end{equation*}
\end{thm}
\begin{pf}
We have with H\"older's inequality:
\begin{equation}
\begin{split}
|\langle(\Delta_h - \tilde{\Delta}_h)\omega,\varphi\rangle| =
\bigl|\sum_{i\in\mathcal{I}}\sum_{j=1}^N\Gamma_if_{ij}\varphi(\vv{x}_j)\bigr| \\
\leq
\lnorm{(\Gamma_i)_{i\in\mathcal{I}}}
\bigl\Vert\bigl(\sum_{j=1}^{N}f_{ij}\varphi(\vv{x}_j)\bigr)_{i\in\mathcal{I}}\bigr\Vert_{l^\infty}.
\end{split}
\end{equation}
Applying the triangle-inequality, Lemma~\ref{lem:fnorm} and the Sobolev embedding
yields the result.\qed
\end{pf}
For the additional error to be of the same order as the error of the full scheme,
we introduce a new user defined constant $C_{\mathrm{diff}}$ and require:
\begin{equation}\label{eqn:reduced_bound}
\lnorm{(\Gamma_i)_{i\in\mathcal{I}}} \leq
C_{\mathrm{diff}}h^{n+2}\lnorm{(\Gamma_i)_i}.
\end{equation}
To minimise the number of diffused particles, we ignore those with the smallest
amount of individual circulation, until this bound is reached. Using the same
methods as above, it is easily verified that the reduced operator does conserve
circulation and linear momentum, however, it does \emph{not} conserve angular
momentum. 

\subsection{Limitations of the Method}\label{subsec:limitations}
After having investigated consistency, stability, and conservation properties
of the method we want to make some comments on its limitations. First, we want
to point out that unlike claimed by Shankar and van Dommelen~\cite{shankar1996}
the vorticity redistribution method \emph{does not} extend to arbitrary orders
of accuracy. The stability proof relies on the fact that the stencils are non-%
negative. While other stable stencils might and probably do exist, we are not
aware of any stability proof. As Seibold~\cite{seibold2006,seibold2008} points
out, any third or higher order method needs to fulfil the moment conditions for
$|\alpha| = 4$. A simple linear combination of these equations yields:
\begin{equation}
\sum_{j\neq i}f_{ij}|\vv{r}_{ij}|^4 = 0,
\end{equation}
which due to the non-negativity constraint can only be fulfilled for
$f_{ij} \equiv 0$. The zero stencil, however, is inconsistent with the moment
conditions for $|\alpha| = 2$. We may point out, however, that the method of
particle strength exchange (PSE) similarly requires a positive kernel function
for its stability proof, thereby equally limiting it to second order
accuracy~\cite{cottet2000}.

Secondly, we point out that the matrix $\sm{F}$ discontinuously depends on
the particle positions: as they move around, they may enter and leave each
other's neighbourhoods, allowing for jumps between zero and non-zero in the
corresponding entries $f_{ij}$. In fact, in general, the solution to the
moment equations is not even unique. It is thus hard to analyse the effect
of higher-order time-stepping schemes on the method's accuracy.

\section{Implementation}\label{sec:implementation}
It has been claimed that the VRM is a slow algorithm, especially when compared
to the PSE scheme, e.\,g., by Cottet and Koumoutsakos~\cite{cottet2000}. On the
other hand, Schrader et al.~\cite{schrader2010} report that their DC-PSE method
also takes up as much as 90\% of total CPU time, and compare its computational
speed with that of the VRM. We believe that the computational cost associated
with the VRM has been greatly overestimated; mostly due to implementation issues.
In this section we discuss some of these issues and illustrate a heuristic which
can further speed up the method significantly. In our final implementation the
velocity computation took about three times longer than the evaluation of the
discrete Laplacian.

\subsection{Solution of the Moment Equations}\label{sec:simplex}
Equation~\eqref{eqn:momcondmv} is a classical \enquote{phase I problem} of the
Simplex algorithm for linear programming problems. When we use an insertion
scheme such as the one described in section~\ref{sec:insertion} this system is
underdetermined, with a fixed, small number of rows $m = 5$ or $m = 9$,
corresponding to the number of moment conditions, and a variable number of
columns, corresponding to the size of the neighbourhood $\neighi\setminus
\lbrace i\rbrace$. 

The theory of simplex algorithms is to vast to be treated in detail here, such
that we can only give a some key remarks and refer to the literature, e.\,g.,
Fletcher's book~\cite{fletcher2000}, for further details. Assuming that the
moment conditions do have a solution, phase I of the algorithm always returns
one with $m$ non-zero entries corresponding to a certain subset of
particles in the neighbourhood. These non-zero entries are called basic
variables. Setting the fractions $f_{ij}$ for the remaining particles to zero,
the solution can be obtained by solving an $m\times m$ linear system. The
simplex algorithm is a systematic, iterative way of finding a valid set of
basic variables. In every iteration of the algorithm an $m\times m$ system
consisting of varying sets of columns of $\sm{V}_i$ needs to be solved,
typically by means of an LU-decomposition, which is of $\bigO{m^3}$
complexity~\cite{golub2013}.

The efficiency of the method thus crucially depends on the number of rows of
$\sm{V}_i$, which should be kept as small as possible. Shankar and van
Dommelen~\cite{shankar1996} use a different linear programming problem,
aiming to minimise the maximum norm of the solution. By doing so, they solve a
problem involving $4m$ rows, effectively making each iteration 64 times more
expensive. One should thus keep the original formulation~\eqref{eqn:momcondmv}.
Furthermore, optimising the solution with respect to some target value forces
to algorithm to enter phase II, which further increases its cost without
improving the method's order of convergence. One might try to optimise the
error constant by choosing an optimisation criterion that favours close
particles. However, in regard of the later introduced heuristic of small
neighbourhoods in section~\ref{sec:small}, it is not immediately clear if this
additional optimisation step is cheaper than choosing smaller values of $h$.

Note that the two possible values of $m$ are very small and fixed. An efficient
implementation should thus make use of this fact: all loops of the LU-decomposition
can be unrolled, enabling compilers to perform aggressive optimisations. The LAPACK
routines, on the other hand, were optimised for larger problems with dynamic,
varying sizes~\cite{lapack}.

There are several approaches to avoid a from-scratch computation of the LU-%
decomposition in every iteration of the method. Updating LU-decompositions
instead of recomputing them, however, typically is only effective for larger
values of $m$: the Fletcher--Matthews update, for example, is reported to be
effective for $m>10$~\cite{fletcher1984}.

Note that the matrix $\sm{V}_i$ is fully populated and---as the number of
neighbours is typically limited---of small to moderate size. On the other
hand, most available implementations of the Simplex algorithm as well as a
substantial part of the available literature focus on large-scale, sparse
problems. In other words, they are optimised for the opposite case and
thus cannot deliver good performance for our problem. Implementing an efficient,
dense simplex method is essential for the overall performance of the VRM. As
this task is not straight forward, some authors, e.\,g., Lakkis and Ghoniem~%
\cite{lakkis2009}, prefer to solve the non-negative least-squares problem
instead:
\begin{equation}
\min_{\sv{f}_i\geq 0}\vert\sm{V}_i\sv{f}_i-\sv{b}_i\vert^2,
\end{equation}
where $\vert\cdot\vert$ refers to the Euclidean norm. This problem can be solved
using the algorithm due to Lawson and Hanson~\cite{lawson1995}, which solves
an \emph{unconstrained} least-squares problem in each iteration. However, the
size of this unconstrained problem varies in every iteration, making it harder
to unroll loops a priori. Additionally, these problems are typically solved using
QR- or LQ-decompositions, which are more expensive than the LU-decomposition. We
therefore do not further pursue this approach.

\subsection{Insertion of New Particles}\label{sec:insertion}
In order to ensure that non-negative stencils exist, particles need to have
sufficiently many neighbours which also need to fulfil certain geometric
conditions. Seibold~\cite{seibold2006,seibold2008} gives the exact conditions
for the first order case $n=1$ as well as the following sufficient condition:
seen from the centre of the neighbourhood, the angle between two adjacent
particles may be no more than 45\textdegree. Assuming a given maximum hole-size
in the particle cloud, he also gives a sufficient upper bound $Rh$ for the
neighbourhood size. These conditions could in principle be implemented in a VRM
scheme, resulting in a strong guarantee that positive stencils always exist.

\begin{figure}
\centering
\includegraphics[scale=0.4]{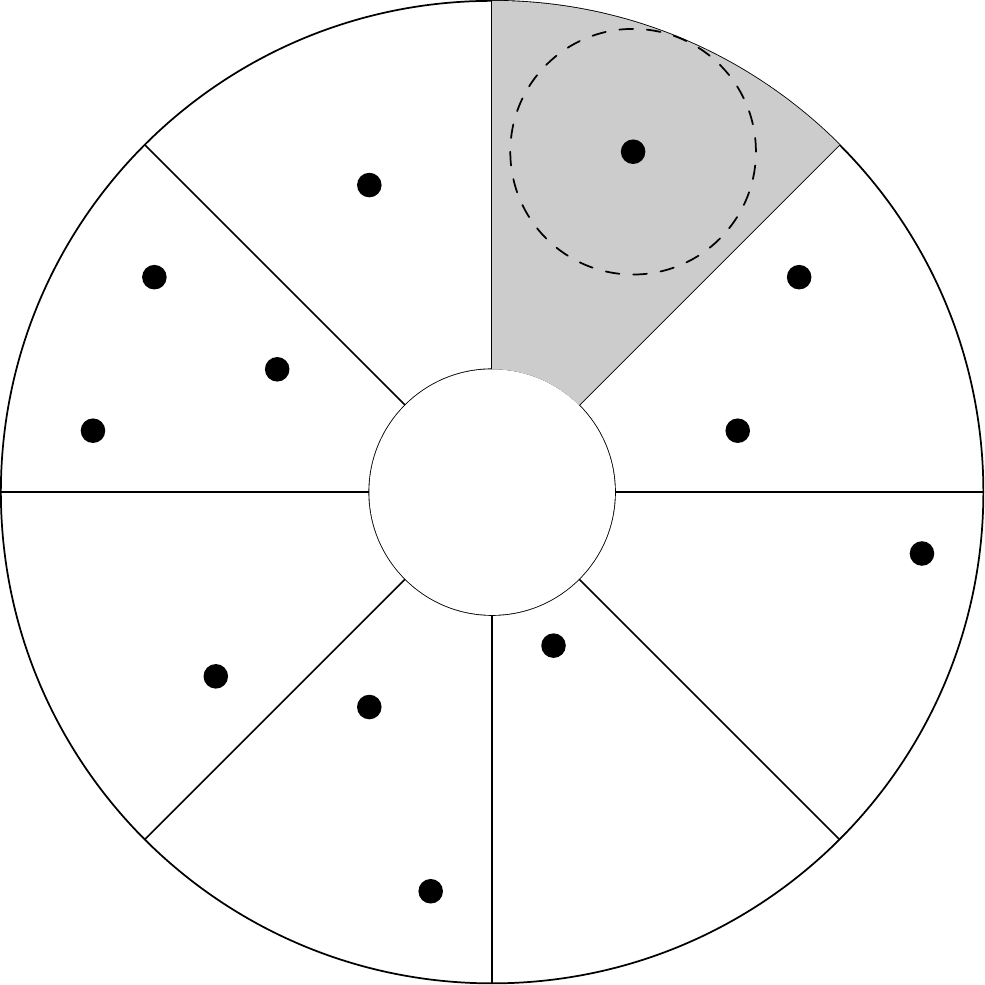}
\caption{\label{fig:insertion}Illustration of a particle neighbourhood and the
insertion strategy. Each of the eight segments except for the shaded one
contained at least one particle. In the shaded segment a new particle is inserted
on the centre line at radial position $1.5h$. No particle can be closer than
$0.5h$ to the newly inserted particle: the circle of that radius is indicated
using a dashed line and is completely included in the previously empty segment.}
\end{figure}

However, as he points out, these conditions are often too strict. We thus pursue
a different approach. Instead of directly checking the angles between each pair
of adjacent particles, we subdivide the neighbourhood into eight segments of
45\textdegree\ each, as illustrated in figure~\ref{fig:insertion}. In order to
avoid wasting computational resources, we do not want to insert new particles
that would violate the lower bound in~\eqref{eqn:neighbourhood} for any other
particle. However, we also want to avoid small values of $r$, to prevent the
time-step constraint~\eqref{eqn:apriobound} from becoming too strict. As a
compromise we choose $r=\nicefrac{1}{2}$ and $R=2$ and apply the following
insertion strategy: if any neighbourhood segment contains no particles, a new
particle is inserted on the segment's centre line at radial position $1.5h$. As
illustrated in figure~\ref{fig:insertion}, this ensures that the newly inserted
particle does not violate any other particle's lower bound on its neighbourhood.

This insertion strategy ensures that particles are at most spaced $2h$ apart.
According to theorem 6.11 of Seibold's thesis, choosing the upper bound of the
neighbourhood size as $R\geq 5.23$ then guarantees the existence of positive
stencils. However, in our numerical experiments, such a large choice was not
necessary and all computations worked well with $R=2$. Experiments conducted
with a slightly rotated reference frame indicated that the results of this
strategy do not significantly depend on the coordinate system used.

Unlike claimed by Cottet and Koumoutsakos~\cite{cottet2000}, insertion of empty
particles is different from remeshing: it leaves the vorticity field \eqref{eqn:vfield}
unchanged, thereby introducing no error and it does not rearrange existing particles.
For this reason the VRM is a truly mesh-free method.

\subsection{Small Neighbourhoods}\label{sec:small}
As pointed out in section~\ref{sec:simplex}, the simplex method systematically
determines a subset of particles leading to a non-negative solution of the
moment equations. One can consequently lower iteration counts by reducing the
number of particles in the neighbourhood. In most cases a non-negative solution
exists if there is just one particle in every 45\textdegree-segment of the
neighbourhood. This leads us to the following approach: for every particle
neighbourhood, choose the closest particle of each segment. We call the resulting
subset the \emph{small neighbourhood}. We then apply the simplex method to this
small neighbourhood. By choosing the segments' closest particles, we aim to
locally reduce $R$, thereby minimising the error constant. Only if no non-negative
solution was found, we retry with the complete neighbourhood. In our numerical
examples, depending on $h$, this only happened in a negligible (less than a
hundred) number of cases.

This approach has the advantage that \emph{all} matrices and vectors involved in
the simplex algorithm can be statically allocated, avoiding the overhead of
dynamic memory allocation and further enabling the compiler to unroll more loops.
In our experiments in section~\ref{sec:speed}, the use of these small
neighbourhoods instead of the complete ones lead to a threefold speed-up.

Note that after the assembly of the Vandermonde matrices $\sm{V}_i$, this
approach leads to a set of completely decoupled, small problems of fixed size.
We thus have an \emph{embarrassingly parallel} problem, making it ideally
suited for computations on many-core processors, such as GPUs or the Intel Xeon
Phi.

\section{Numerical Experiments}
As Shankar and van Dommelen point out in their work~\cite{shankar1996}, the
Lamb--Oseen flow is an ideal test-case for vortex particle methods: its initial
condition is a single Dirac delta distribution:
\begin{equation}
\omega(0,\vv{x}) = \Gamma\delta(\vv{x}),
\end{equation}
and can thus be exactly represented in a vortex particle method. The analytic
solution is infinitely smooth and valid for the heat-equation~\eqref{eqn:heat_eq}
as well as the vorticity-transport equation~\eqref{eqn:vte}:
\begin{equation}
\omega(t,\vv{x}) = \frac{\Gamma}{4\pi\nu t}e^{-\frac{|\vv{x}|^2}{4\nu t}}\qquad (t>0).
\end{equation}
The corresponding velocity field is given by:
\begin{equation}
\vv{u}(t,\vv{x}) = \frac{\Gamma}{2\pi\vert\vv{x}\vert}
\biggl(1-\exp\biggl(-\frac{\vert\vv{x}\vert^2}{4\nu t}\biggr)\biggr)\widehat{\gv{\varphi}},
\end{equation}
where $\widehat{\gv{\varphi}}$ refers to the unit vector in circumferential
direction at position $\vv{x}$. In the following, we will describe several
numerical experiments carried out on this flow. Mimicking Shankar and van
Dommelen's case of $\mathrm{Re} = 50$, we chose $n=1$, $\Gamma = 2\pi$,
$C_\mathrm{diff} = 1$, $\nu=\nicefrac{1}{50}$. We choose higher resolutions,
however, and stop time-integration at $t = 1$.

\subsection{Convergence with respect to h}
We consider the cases with and without convection, corresponding to the
Navier--Stokes equation and the heat equation, respectively. In the case of
the heat equation, we use Euler's method to advance the solution in time
and choose a fixed time step:
\begin{equation}\label{eqn:viscstep}
\Delta t = \frac{1}{8} \frac{(rh)^{2}}{4\nu}.
\end{equation}

As mentioned previously, in vortex methods it is customary to replace the singular
Biot--Savart kernel $\vv{K}$ with a regularised one $\vv{K}_\varepsilon$. We use
the following second order kernel obtained after Gaussian smoothing:
\begin{equation*}\label{eqn:smoothkernel}
\vv{K}_\varepsilon(\vv{x}) = \frac{(-y,x)^\top}{2\pi\vert\vv{x}\vert^2}
\biggl(1-\exp\biggl(-\biggl\vert\frac{\vv{x}}{\varepsilon}\biggr\vert^2\biggr)\biggr),
\qquad \vv{x} = (x,y)^\top.
\end{equation*}
Our particle insertion strategy guarantees that particles are at most spaced
$2h$ apart. To ensure sufficient overlap we choose $\varepsilon = 3h$.
A fast multipole method (FMM) similar to that of Dehnen~\cite{dehnen2002} of
order $p = 16$ and multipole acceptance criterion $\theta\leq0.8$ is used
to speed up the velocity computation.

Practical experience has shown, that higher order time-stepping methods are 
required to maintain linear and angular momentum in the case of enabled convection.
Like Shankar and van Dommelen, we choose the classical Runge--Kutta method (RK4)
in this case. In order to resolve particle movement accurately, the time-step
is adaptively chosen as the minimum of $\eqref{eqn:viscstep}$ and the following
CFL-type condition:
\begin{equation}
\Delta t\leq \frac{1}{8}\min_{i=1,\ldots,N}{\frac{h}{|\vv{u}_i|}}.
\end{equation}
We want to stress that this second bound is not required to ensure stability:
experiments without this restriction showed no instabilities and gave reasonable
results, however, the errors in linear and angular momentum were larger.

As it is difficult to compute Sobolev-norm $\qnorm{\cdot}$ explicitly, we try
to approximate the $L^2$-error of the corresponding velocity. As the system
contains infinite energy, we need to limit the area of integration. We chose
$A = [-1.5,1.5]^2$, as all particles were contained within this region. By
means of numerical quadrature we then evaluate:
\begin{equation}
e_{\vv{u}} = \frac{\Vert\vv{u} - \vv{u}_h\Vert_{L^2(A)}}{\Vert\vv{u}\Vert_{L^2(A)}},
\end{equation}
where $\vv{u}_h$ stands for the velocity field which is obtained from the particle
approximation for a chosen value of $h$ using the smoothed kernel $\vv{K}_\varepsilon$.

\begin{figure}
\centering
\includegraphics[width=\textwidth]{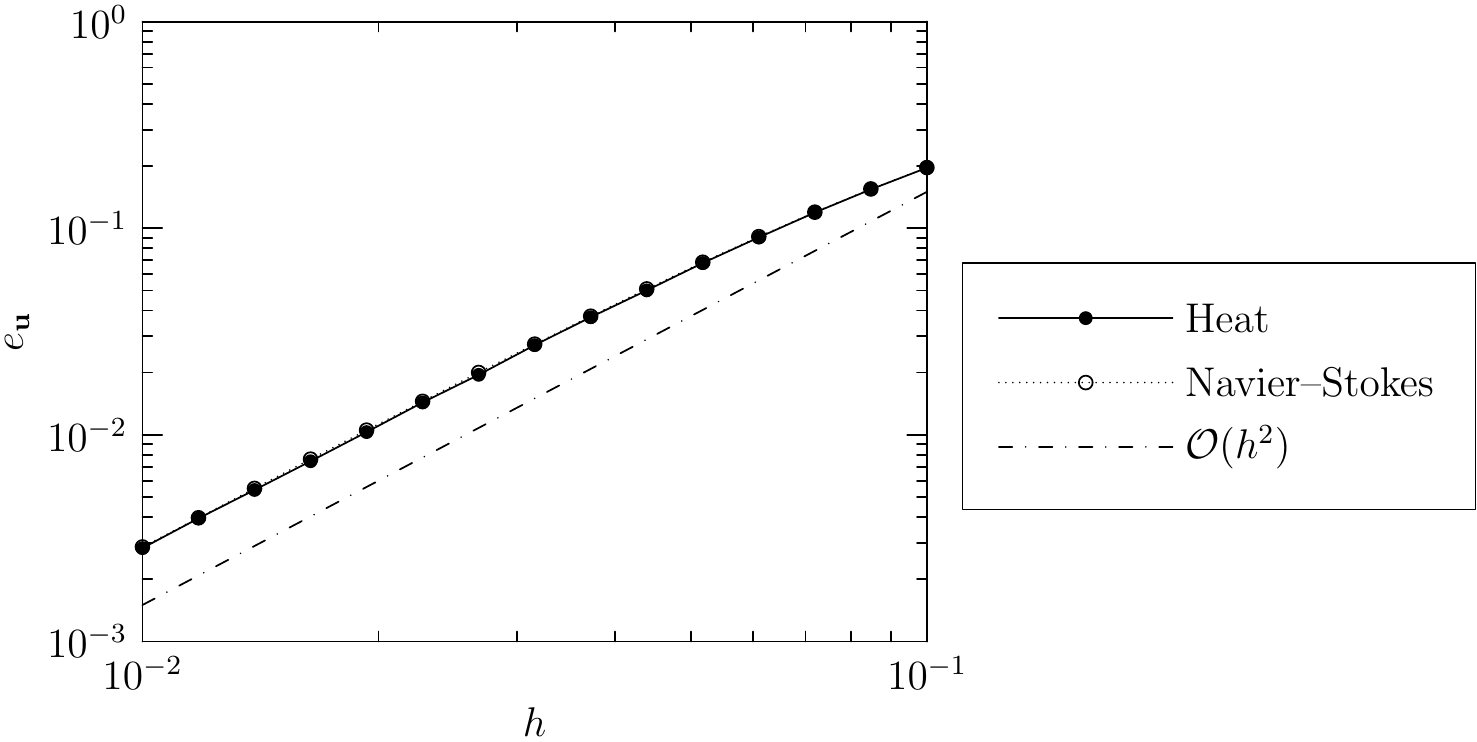}
\caption{\label{fig:errorplot}Error estimates for the heat and Navier--Stokes
equations for varying values of $h$ at $t=1$. Their values essentially coincide
and exhibit an $\bigO{h^2}$ convergence behaviour.}
\end{figure}

Figure~\ref{fig:errorplot} shows the observed error estimates for various
values of $h$. Even though the expected convergence rate was $n=1$, we actually
observe second order convergence behaviour. This is similar to the observations
by Seibold, who explains this using a symmetry argument: the classical five-%
point finite-difference stencil achieves second order accuracy due to the
symmetry of the particle locations. However, the insertion strategy and the
definition of the particle neighbourhoods preclude extreme cases of asymmetry,
which might result in the observed second order convergence. Seibold, however,
does not exclude particles according to equation~\eqref{eqn:reduced_bound}.
It is thus comes as a surprise that even the reduced operator exhibits this
behaviour. As both curves form a nearly straight line and essentially coincide,
we suspect that the smoothing error dominates for this choice of parameters.

\begin{figure}
\centering
\includegraphics[width=\textwidth]{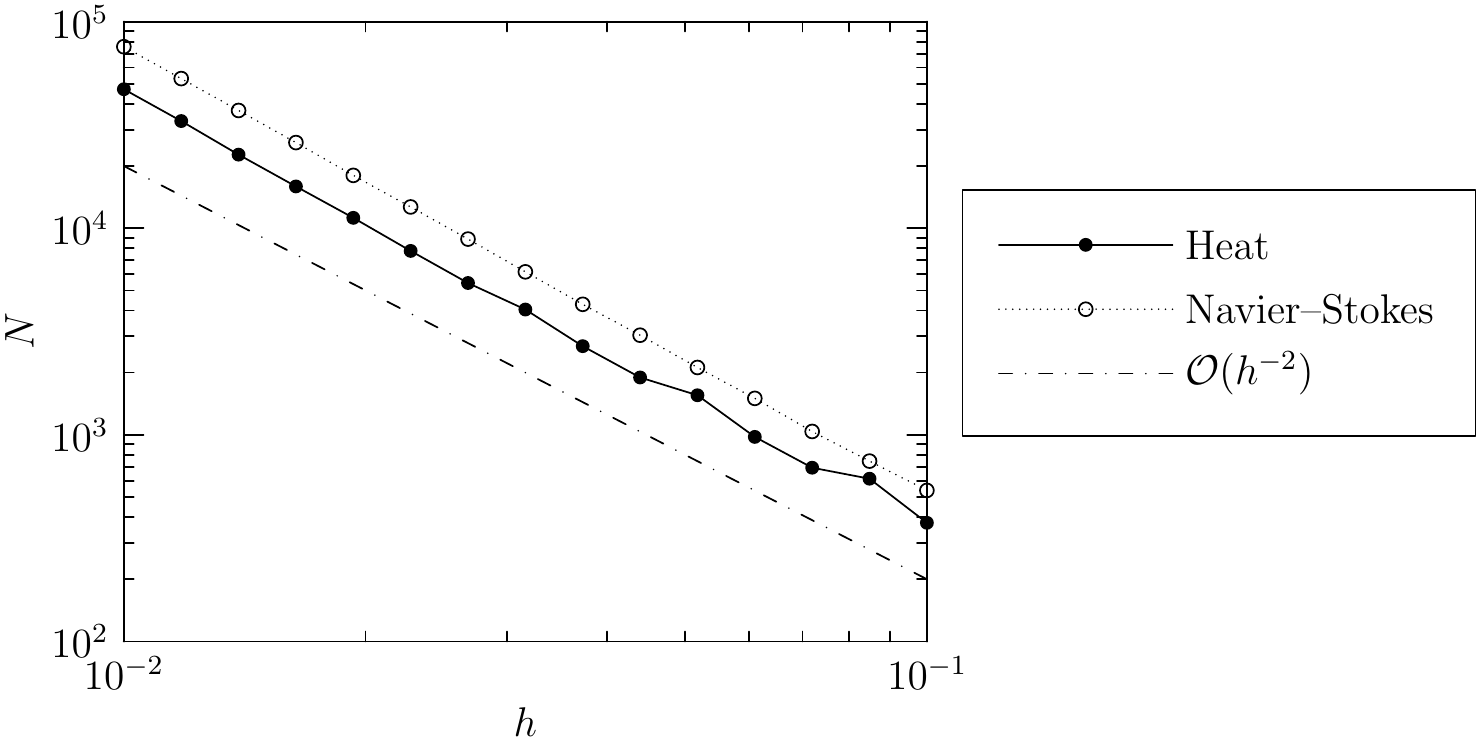}
\caption{\label{fig:particles}The number of particles in the final step of the
computation for the heat and Navier--Stokes equations. The curves show a particle
growth that scales as $\bigO{h^{-2}}$, despite the fact that equation~%
\eqref{eqn:reduced_bound} is getting stricter for decreasing mesh-sizes. The
ratio between the two curves' values remains approximately fixed at around
$1.6$.}
\end{figure}

Figure~\ref{fig:particles} shows the number of particles in the final time-step
of the computation. It increases approximately as $\bigO{h^{-2}}$, as one would
expect in a grid-based computation. This again is surprising, as bound~\eqref{%
eqn:reduced_bound} gets stricter for decreasing $h$. Due to the convection in
the Navier--Stokes case, more particles need to be inserted as they move around.
In our simulation, this caused an increase in the number of particles of a nearly
constant factor 1.6. 

\begin{figure}
\centering
\includegraphics[width=\textwidth]{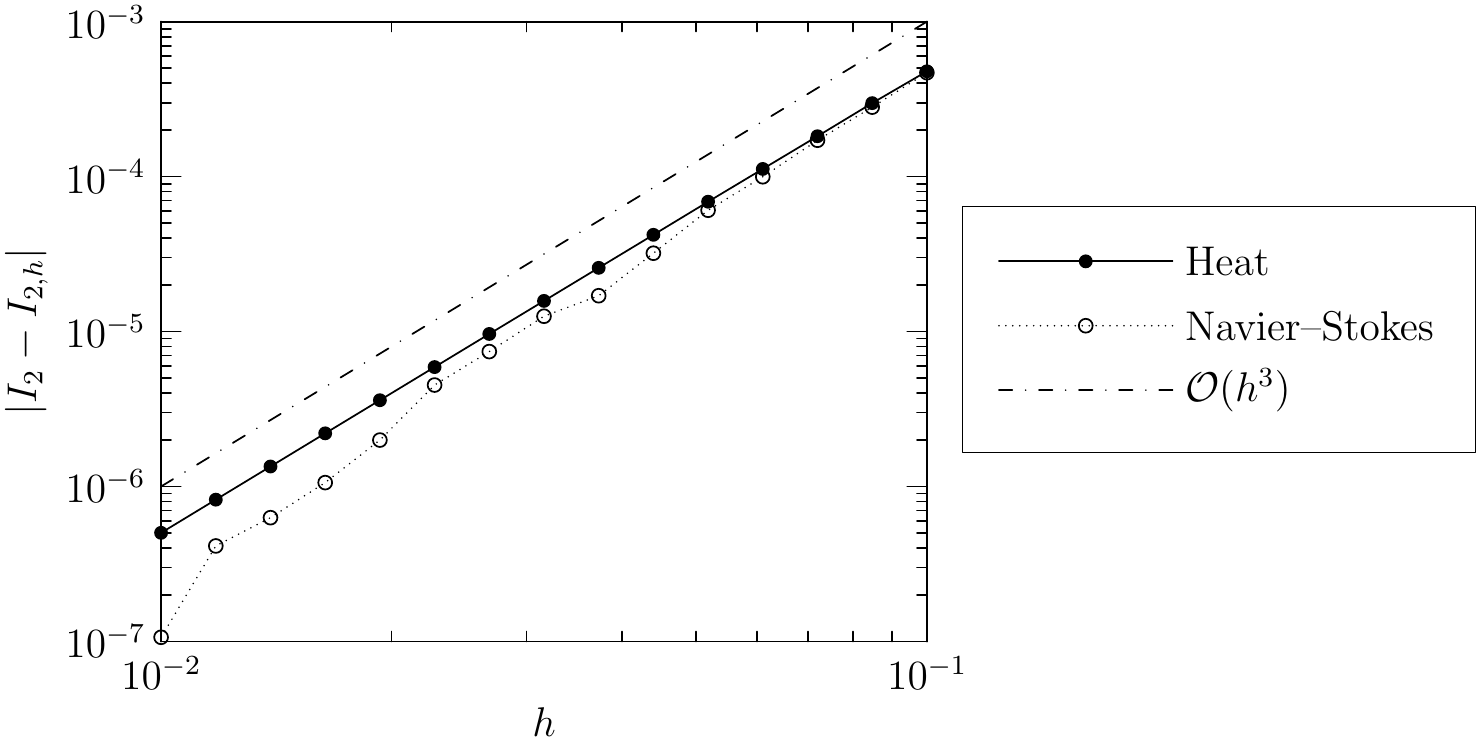}
\caption{\label{fig:angmomerr}Error in angular momentum at the final time-step
for the heat and Navier--Stokes equations. The error decays at a rate of 
$\bigO{h^3}$, the same exponent as in condition~\eqref{eqn:reduced_bound}. In
case of the Navier--Stokes equations, the error decreases even faster, in a
less clear-cut manner.}
\end{figure}

As shown in section~\ref{subsec:reduced}, the reduced operator conserves circulation
and linear momentum exactly. In the case of the heat equation this remains the
case when a time-stepping scheme is applied: the error in $I_0$ and $\vv{I}_1$
was of the order of the machine accuracy. For the Navier--Stokes equation this
is only true for the circulation. For all choices of $h$ the error in linear
momentum varied between $\bigO{10^{-6}}$ and $\bigO{10^{-7}}$. We believe this to
be a result of the limited accuracy of the FMM code that was used for the velocity
computation and the error introduced by the Runge--Kutta method.
Figure~\ref{fig:angmomerr} shows the error in angular momentum $I_2$. The values for
the heat equation decrease at a rate of $\bigO{h^3}$, similar to the bound~%
\eqref{eqn:reduced_bound}. In the convective case the error decays somewhat faster,
in a less clear-cut manner. We believe this to be a result of the increased number
of particles. We thus conclude that for the chosen values of $h$, the error in
angular momentum induced by using the reduced operator~\eqref{eqn:reducedop}
dominates that of the FMM and the time-stepping scheme.

\begin{figure}
\centering
\includegraphics[width=.5\textwidth]{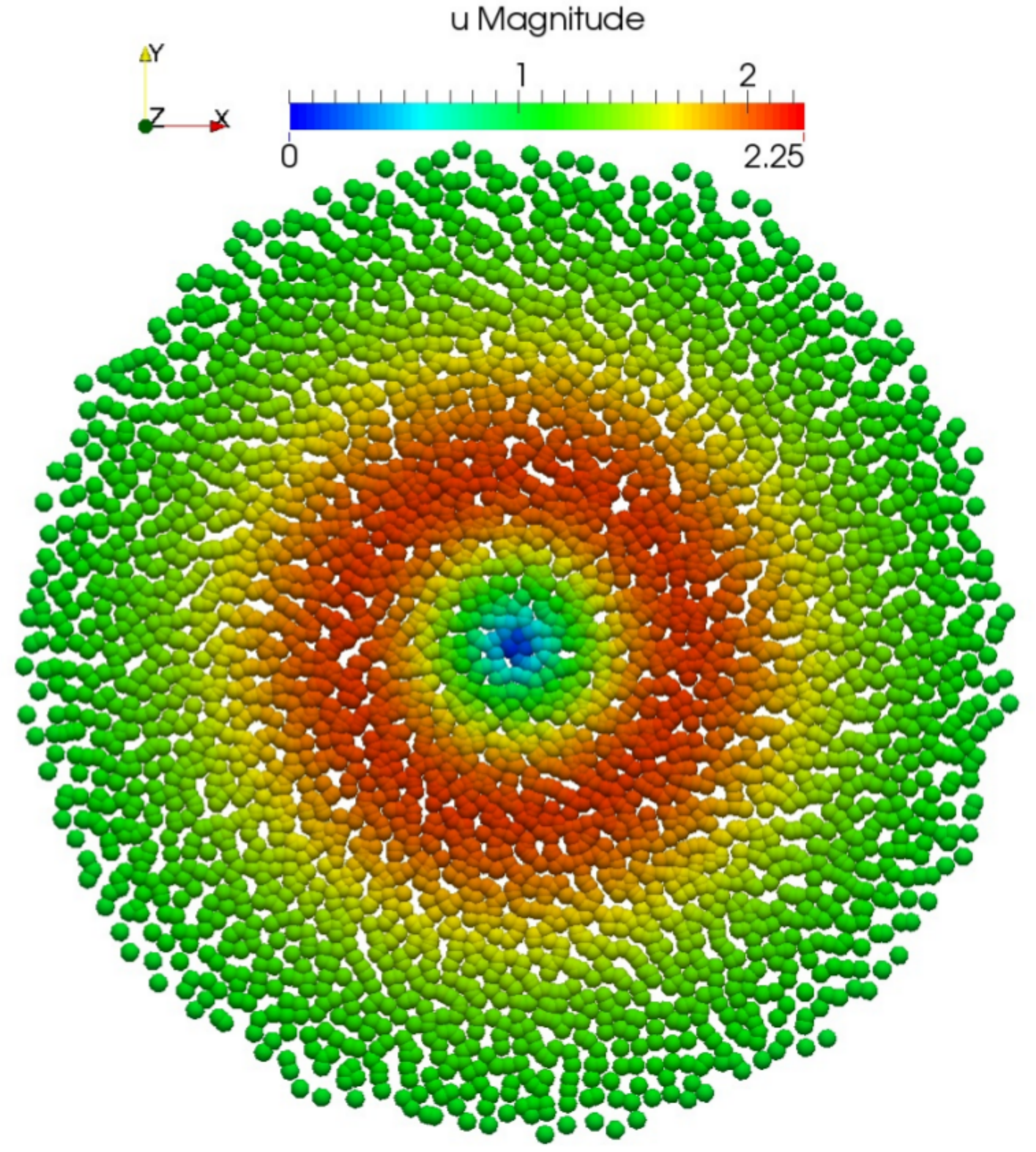}
\caption{\label{fig:solution}Plot of the smoothed velocity at the particle
locations for $h=0.04$ at $t=1$. Despite the asymmetric particle distribution,
caused by the convection, the velocity field remains very symmetric. The particle
cloud takes the shape of a circle. Even in this under-resolved case, the method
yields qualitatively good results.}
\end{figure}

Figure~\ref{fig:solution} shows the velocity at the particle locations for $h=0.04$
at $t=1$ with enabled convection. Despite the asymmetry in the particle locations
caused by the convection, one can see that the velocity field remains quite symmetric.
The reduced operator prevents the creation of particles that would carry insignificant
amount of circulation. For this reason, the particle cloud takes the shape of a circle
around the origin: vorticity decays exponentially with the distance to the origin. At
$\varepsilon = 3h = 0.12$ the resolution is not high enough to accurately represent
steep velocity gradient at the centre of the flow. However, due to the good
conservation properties, we obtain a qualitatively good solution already at this
under-resolved computation.

\subsection{Computational Speed}\label{sec:speed}
In order to assess the speed of the method, we measured the time needed
to evaluate the velocity and the Laplacian for $h=0.01$. For the Laplacian,
we compared the performance of two codes: the first code uses LAPACK to
decompose the arising linear systems in each simplex iteration and takes the
complete particle neighbourhood into account. The second code uses small
neighbourhoods as described in section~\ref{sec:small} and an implementation
using completely unrolled loops in the LU-decomposition. The code was
parallelised using OpenMP, where task-based parallelism was used for the FMM.

\begin{figure}
\centering
\includegraphics[width=\textwidth]{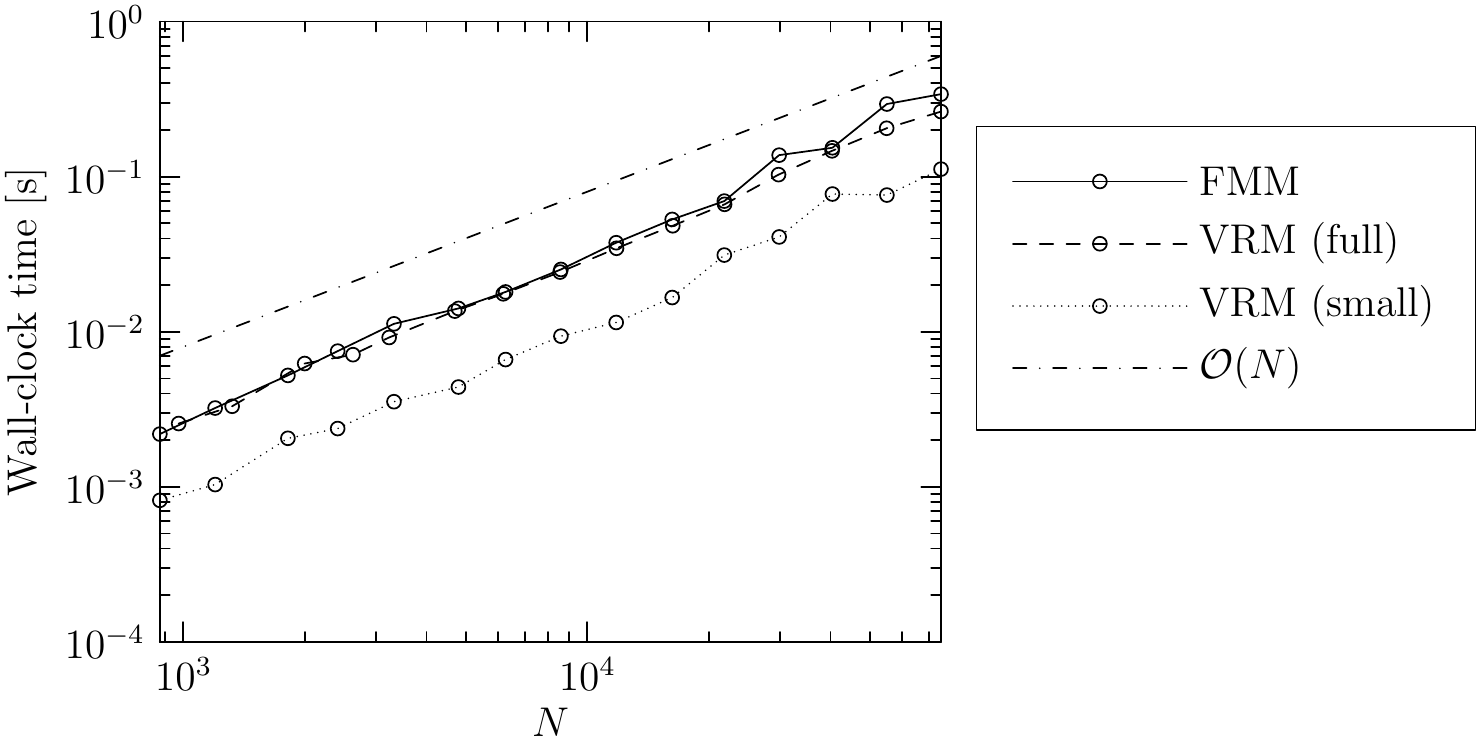}
\caption{\label{fig:speed}Required CPU time for the VRM with the full and
small neighbourhoods in comparison to the FMM. The computations were performed
on an Intel Xeon E5-1650v3, a six-core processor running at 3.5\,GHz. The line
corresponding to the FMM is jagged due to the task-based parallelism used in
the implementation. The VRM computation can be greatly accelerated using small
neighbourhoods, it is then about three times faster than the corresponding velocity
computation.}
\end{figure}

Figure~\ref{fig:speed} shows the required time for each computation depending on
the number of particles involved. One can see that all computations scale linearly
with $N$, however, with different constant factors. The code using small
neighbourhoods performs about three times faster than the corresponding code using
the complete ones. This clearly highlights the benefit of trying small
neighbourhoods first. It also performs about three times as fast as the
corresponding FMM code. Further measurements showed that, in the case of small
neighbourhoods, only about one third of the time was used for the actual simplex
solver, while the remaining time was spent finding neighbourhoods and inserting
new particles. A hash-based algorithm was used for this, causing the resulting
curve to be jagged due to caching effects.

Note that these numbers cannot be directly compared to those reported by Shankar
and van Dommelen: they compare a single VRM computation to that of a convective
step performed using the Runge--Kutta method, i.\,e., involving four velocity
computations. In this setting, their VRM computation takes about 25\% longer than
the convective step, i.\,e., five times longer than a single velocity evaluation.
In comparison to the respective FMM codes, our VRM computation thus is about 15
times faster.

\section{Conclusion and Outlook}\label{sec:conclusion}
\subsection{Conclusion}
We have introduced a splitting-free variant of the vorticity redistribution
method (VRM). Using the new concept of small-neighbourhoods, its speed compared
to the original method can be greatly accelerated and typically is below that
of the corresponding velocity computation. Equation~\eqref{eqn:reduced_bound}
allows us to efficiently and consistently reduce the number of diffused particles.
We have illustrated that the method can be implemented efficiently and that
previous claims on the slow speed of the VRM are probably due to implementation
issues. The large number of small, independent, fixed-size problems involved
makes it an ideal candidate for parallelisation on coprocessors such as GPUs or
the Intel Xeon Phi. We conclude this text with a few possible extensions on the
method.

\subsection{Outlook}\label{sec:outlook}
In light of the quadratic time-step bound~\eqref{eqn:apriobound}, an interesting
topic for future research might be the application of implicit time-stepping
schemes in periodic flows. As the convective part of the equations is non-stiff,
this seems to be an ideal use-case for IMEX multistep schemes~\cite{ascher1993}.
After having convected the particles, $\sm{F}$ could then be readily assembled,
leading to a linear system. As Seibold discusses in his work~\cite{seibold2010},
due to the positivity and sparsity of the stencils, such systems can effectively
be solved using algebraic multigrid methods.

The definition of a particle's neighbourhood in equation~\eqref{eqn:neighbourhood}
excludes particles that are too close to that particle. In order to save
computational resources, it may thus be desirable to remove particles in areas
where they get too close to one another. Instead of approximating the Laplacian
as described in this article, one can apply the same methodology to approximate
the identity operator using a particle's neighbours. This way, a particle can be
redistributed to its neighbours and subsequently be removed. Lakkis and Ghoniem
\cite{lakkis2009} successfully applied a similar procedure and reported a
significant reduction in the number of particles. 

Finally, we would like to conclude this text by thanking the editor and reviewers
for their comments which helped improving the quality of this article. This work
was financially supported by the Keio Leading-edge Laboratory of Science and 
Technology (KLL). The first author also receives the MEXT scholarship of the
Japanese Ministry of Education.

\section*{References}
\bibliographystyle{elsarticle-num}
\bibliography{literature}
\end{document}